\documentclass{amsart}

\usepackage{amsfonts}
\usepackage{amsthm}
\usepackage{amscd}

\usepackage{graphics}
\usepackage{epsfig}		%
\usepackage{graphicx}		%
\usepackage{psfrag}

\theoremstyle{definition}

\theoremstyle{definition}

\theoremstyle{definition}

\newcommand{\Z}{{\mathbb Z}}

\newcommand{\R}{{\mathbb R}}

\renewcommand{\S}{\Sigma}

\newcommand{\s}{\sigma}

\renewcommand{\r}{\rho}
\newcommand{\D}{\Delta}

\newcommand{\G}{\Gamma}

\renewcommand{\L}{\Lambda}

\newcommand{\iso}{\cong}
\newcommand{\bdry}{\partial}

\renewcommand{\int}{\text{Int}\ }

\def\hpic #1 #2 {\mbox{$\begin{array}[c]{l} \epsfig{file=#1.eps,
height=#2} \end{array}$}}
\def\frac#1#2{{#1\over#2}}

\begin{document}

\title{Configuration Spaces of Colored Graphs\footnote{This article was published in \emph{Geometriae Dedicata}, vol.~92 (2002), pp.~185--194.}}
\author{Aaron Abrams
}
\thanks{Partially supported by the frequent, positive, and humorous influence of John Stallings.}
\dedicatory{This paper is dedicated to John Stallings on his 
64.8$\,^{th}$ birthday.}

\begin{abstract}

This paper is intended to provide concrete examples of
concepts discussed elsewhere in this volume, especially
splittings of groups and non-positively curved cube complexes but also
other things.  The idea
of the construction (configuration spaces) is not new, but 
this family of examples doesn't seem to be well-known.  
Nevertheless they arise in a variety of contexts; applications 
are discussed in the last section.  Most proofs are omitted.

\end{abstract}

\maketitle

\section{Graphs}

\subparagraph{\bf 1.1.}
A \emph{graph} $G$ is a 1-dimensional cell complex; a cell $\s$ of
$G$ is either a \emph{vertex} (0-cell) or an \emph{edge} (1-cell).  
We always
assume $G$ is locally finite, i.e. each vertex has only finitely 
many edges attached to it.  We do not need to distinguish between
$G$ and its underlying topological space; thus we refer both 
to points of $G$ and to cells of $G$.

\subparagraph{\bf 1.2.}
The ``usual'' (topological) $n$-point configuration space of $G$
is the space of $n$-tuples of distinct points of $G$:
\begin{eqnarray*}
C_n^{\rm top}(G) &= &\{(x_1,\ldots,x_n)\in G\times\cdots\times G
\ \big|\ x_i\ne x_j \text{ if } i\ne j\} \\
&= &(G\times\cdots\times G) - \D,
\end{eqnarray*}
where $\D$ is the ``diagonal.''  The fundamental group of 
$C_n^{\rm top}(G)$
probably deserves to be called the \emph{pure braid group} of $G$, but
we reserve this word for something slightly different (see {\bf 4.1}).

Note that $C_n^{\rm top}(G)$ is (usually) noncompact, and it doesn't have
a natural cell structure since the diagonal cuts through many cells of
$G^n$.

\subparagraph{\bf 1.3.}
We define a \emph{``discretized''} or \emph{``combinatorial'' $n$-point 
configuration
space} as follows.  As $G$ is a cell complex, so is the product 
$G^n=G\times\cdots\times G$.  
Let $C_n(G)$ be the largest subcomplex of 
$G^n$ which is contained in $C_n^{\rm top}(G)$; this is also the largest
subcomplex of $G^n$ which doesn't intersect $\D$.  Observe that $C_n(G)$
is a cube complex (since $G^n$ is), and it is compact if $G$ is.
It is a theorem \cite{thesis} that $C_n^{\rm top}(G)$ deformation
retracts to $C_n(G)$ if the edges of $G$ are subdivided enough.

%

\subparagraph{\bf 1.4}
Here is a more positive definition of $C_n(G)$.  
Let $\s$ be a cell (vertex or edge) of $G$; if $\s$ is an edge, define
$\bdry\s$ to be the image of
the attaching map of $\s$ in $G$, and if $\s$ is a vertex, let 
$\bdry\s=\{\s\}$.
Then $C_n(G)$ is the union of those product cells
$\s_1\times\cdots\times\s_n$
of $G^n$ satisfying $\bdry\s_i\cap\bdry\s_j=\emptyset$ for all $i\ne j$.

\subparagraph{\bf 1.5.}
A little notation before we see some examples:
the graph $K_n$ is the 1-skeleton of a simplex of dimension $n-1$,
and the graph $K_{m,n}$ is the join of a set of $m$ points
with a set of $n$ points.

\subparagraph{\bf 1.6.  Example.}
Consider the triangle $K_3$ shown below.  The space $C_2(K_3)$ is one
dimensional, since there do not exist two disjoint edges in $K_3$.
It is connected, with 6 vertices $(A,B)$, $(A,C)$, $(B,A)$,
$(B,C)$, $(C,A)$ and $(C,B)$ and 6 edges $a\times A$, $A\times a$,
$b\times B$, $B\times b$, $c\times C$, and $C\times c$.  The edge
$a\times A$ is
attached to the vertices of $\bdry a\times A$, which are $(B,A)$ and $(C,A)$,
and so on.  Thus $C_2(K_3)$ is the hexagon pictured at right.

\bigskip
\begin{figure}[h] \begin{center}
        \psfrag{a}{\scalebox{2.5}{$a$}}
        \psfrag{b}{\scalebox{2.5}{$b$}}
        \psfrag{c}{\scalebox{2.5}{$c$}}
        \psfrag{A}{\scalebox{2.5}{$A$}}
        \psfrag{B}{\scalebox{2.5}{$B$}}
        \psfrag{C}{\scalebox{2.5}{$C$}}
        \scalebox{.4}{\includegraphics{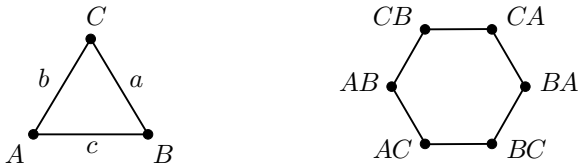}}
        \label{figC2K3}
	\caption{The graph $K_3$ and the space $C_2(K_3)$.}
        \end{center}
\end{figure}

\subparagraph{\bf 1.7.}  One might think of $C_n(G)$ as follows.  Suppose
$n$ (labelled) tokens sit on distinct vertices of $G$, and each 
token may move (along an edge) from
its vertex to an adjacent unoccupied vertex.  A collection of moves can be
executed simultaneously only if they can be done in any order.  With these
rules, the set of possible configurations of the tokens is $C_n(G)$.

\begin{figure}[h] \begin{center}
        \psfrag{A}{\scalebox{2}{$A$}}
        \psfrag{B}{\scalebox{2}{$B$}}
        \psfrag{C}{\scalebox{2}{$C$}}
        \psfrag{O}{\scalebox{2}{$O$}}
        \scalebox{.5}{\includegraphics{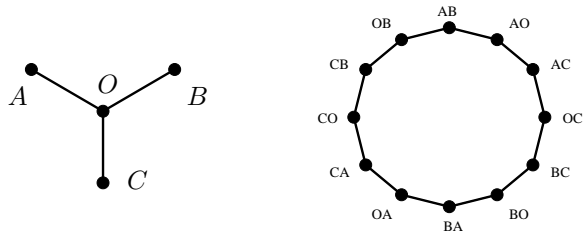}}
	\caption{The graph $K_{1,3}$ and the space $C_2(K_{1,3})$.}
        \label{figC2K13}
        \end{center}
\end{figure}

\subparagraph{\bf 1.8.  Example.}
\label{exC2K13}  The space $C_2(K_{1,3})$ also has no
2-cells.  It has 12 vertices, each of degree 2, and it is connected;
thus $C_2(K_{1,3})$ is a 12-gon.
Exercise:  Place two tokens on $K_{1,3}$ and move them so as to trace
out the loop which is $C_2(K_{1,3})$.  See Figure \ref{figC2K13}.

\bigskip

\begin{figure}[h]\begin{center}
\psfrag{a}{\scalebox{1}{(a)}}
\psfrag{b}{\scalebox{1}{(b)}}
\scalebox{.8}{\includegraphics{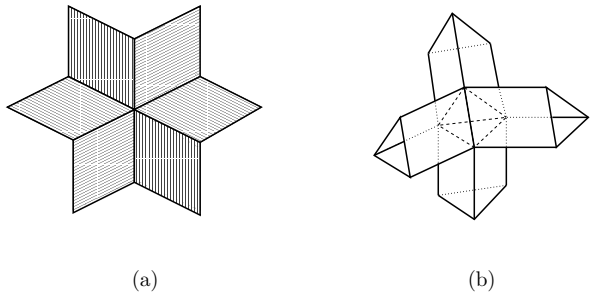}}
\caption{(a)  A neighborhood of a vertex of $C_2(K_5)$.  (b)  The dotted
$K_4$ is $p_1^{-1}(A)$; the rest is a neighborhood in $C_2(K_5)$.  Notice 
that each vertex looks like Figure (a).}
\label{C2K5.fig}
\end{center}\end{figure}

\subparagraph{\bf 1.9.  Example.}
Things get more interesting when we consider
$G=K_5$, with vertices $A$, $B$, $C$, $D$, $E$.
The space $C_2(K_5)$ has $5\cdot 4=20$ vertices.  An arbitrary
vertex, say $(A,B)$, has degree six; it is adjacent to $(A,C)$, 
$(A,D)$, $(A,E)$,
$(C,B)$, $(D,B)$, and $(E,B)$.  The vertex $(A,B)$ is also part of exactly
six 2-cells, which the reader should verify fit together as in Figure
\ref{C2K5.fig}(a).

Consequently, an edge (say $(A,BC)$) of $C_2(G)$ touches exactly two 2-cells
(namely $(AD,BC)$ and $(AE,BC)$).  Thus $C_2(K_5)$ is homeomorphic to
a closed 2-manifold $\S$.  
The excess angle at each vertex is $\pi$, so
the Euler characteristic of $\S$ is $-10$.  (This is just as easily
computed by counting cells.)  It is easy to check that $\S$ is orientable;
hence $\S$ has genus 6.

Let $p_1:C_2(K_5)\to K_5$ be the projection 
onto the first coordinate.  Then $p_1^{-1}(A)$ is a copy of $K_4$ 
(with vertices $(A,B)$, $(A,C)$, $(A,D)$, $(A,E)$), and
a neighborhood of
$p_1^{-1}(A)$ looks
like Figure \ref{C2K5.fig}(b). 
Hence, amazingly, if $K_5$ is embedded in $\R^3$ then
$C_2(K_5)$ can be viewed as the boundary of a regular neighborhood of
the graph (with a suitable cell structure).  

This example, and the next, were discovered at least as early as 1970
\cite{copelandandpatty}.


\subparagraph{\bf 1.10.  Example.}
The reader will enjoy verifying that $C_2(K_{3,3})$ is also 
homeomorphic to a closed orientable 2-manifold.  This time 
some vertices have degree 4 and are ``flat,'' 
whereas others have degree 6 as in Figure \ref{C2K5.fig}(a).  
The manifold $C_2(K_{3,3})$ has genus 4, and it can
be realized as the boundary of a regular neighborhood of the graph
$K_{3,3}$ embedded in $\R^3$.  See also Section \ref{more}.

\section{Colored Graphs}

\subparagraph{\bf 2.1.}
One unfortunate property of configuration spaces of graphs is that
they don't behave well under morphisms $f:G\to L$.  (By a \emph{morphism}
we mean a combinatorial map:  in addition to $f$ being continuous, 
vertices go to vertices, and each edge of $G$
maps homeomorphically to a single edge of $L$.)  For instance, if
$f$ is not injective, then there are vertices $v$, $w$ of $G$ with
$f(v)=f(w)$, and although $(v,w)\in C_2(G)$, there is no sensible
image of $(v,w)$ in $C_2(L)$.

\subparagraph{\bf 2.2.}
We shall overcome this misfortune by introducing colorings.
Let $\mathcal C$ be a (universal) set of colors.
A {\em colored graph} $\Gamma$ is a pair $(G,\varphi)$ where $G$ is
a graph and $\varphi$ is a function from the set of 
vertices of $G$ to $\mathcal C$.  

Following {\bf 1.4}, we define the \emph{($n$-point) configuration space} 
$C_n(\G)$ of the colored graph
$\G=(G,\varphi)$ to be the union of those product cells $\s_1
\times\cdots\times\s_n$ of $G^n$ satisfying 
$$\varphi(\bdry\s_i)\cap\varphi(\bdry\s_j)=\emptyset\quad
\mbox{ for all }\ i\ne j.$$
We call such cells \emph{color-disjoint}.

If the coloring $\varphi$ is injective then $C_n(\G)=C_n(G)$.

\begin{figure}[h]\begin{center}
\scalebox{.5}{\includegraphics[angle=60]{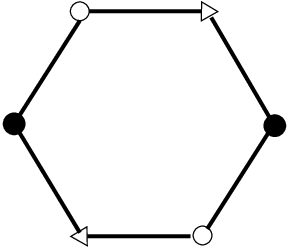}}
\caption{A colored graph $\G$.  The space $C_2(\G)$ is a the disjoint
union of two 12-gons.}
\label{coloredhex}
\end{center}
\end{figure}

\subparagraph{\bf 2.3.  Example.}
Let $G$ be a cycle of length 6, colored as shown in Figure
\ref{coloredhex} with three different
colors.  Note that $C_2(\G)$ contains no 2-cells (unlike $C_2(G)$).
In fact, as the reader should verify, $C_2(\G)$ is a
disjoint union of
two cycles of length 12; it is a (disconnected)
4-fold cover of the space $C_2(K_3)$ from {\bf 1.6}.  See also
{\bf 2.5}.


\subparagraph{\bf 2.4.}
Now let $\G=(G,\varphi_G)$ and $\L=(L,\varphi_L)$ be colored graphs.
A \emph{map of colored graphs} $f:(G,\varphi_G)\to(L,\varphi_L)$ is a 
morphism ({\bf 2.1}), also called $f$, from the graph $G$ to
the graph $L$ such that $\varphi_G=\varphi_L\circ f$.  Such a
map $f$ does induce maps on the colored configuration spaces.  
It is an exercise to prove the following theorem.

\subparagraph{\bf 2.5.  Theorem.}
Let $f:\G\to\L$ be a map of colored graphs.  Then
$f$ induces maps $f^n:C_n(\G)\to C_n(\L)$ for all $n$
by $f^n(x_1,\ldots,x_n)=(fx_1,\ldots,fx_n)$.
Moreover $f^n$ is cellular, and:
\begin{itemize}
\item[(i)]  If $f$ is injective then so is $f^n$;
\item[(ii)]  If $f$ is an immersion then so is $f^n$;
\item[(iii)]  If $f$ is a covering then so is $f^n$;
\item[(iv)]  If $f$ is surjective then so is $f^n$.
\end{itemize}
If $f$ is a $k$-fold cover, then $f^n$ is a $k^n$-fold cover, and
if $f$ is a regular cover with deck transformation group $\pi$, then
$f^n$ is a regular cover with deck transformation group $\pi^n$.

\subparagraph{\bf 2.6.  Example.}
$C_2(\,\mbox{any colored } k\mbox{-fold cover of } K_5)$ is a $k^2$-fold
cover of $C_2(K_5)$; thus it is a closed orientable surface of 
genus $1+5k^2$.  Likewise a $k$-fold colored cover of $K_{3,3}$ gives
a surface of genus $1+3k^2$.

\section{Curvature}

\subparagraph{\bf 3.1}  As the abstract suggests, the spaces
$C_n(\G)$ support metrics of non-positive curvature, in the sense 
that their universal covers satisfy the CAT(0) inequality 
\cite{bridson}.
This is readily proved by checking Gromov's \emph{link condition} 
for cube complexes:

\begin{description}
\item{}A cube complex $X$ has non-positive curvature if and only if for
each cell $\s$ of $X$, every triangle in the link of $\s$ bounds
a 2-simplex in the link of $\s$.
\end{description}

\subparagraph{\bf 3.2.}  We remark that the spaces $C_n(\G)$ are 
generally not negatively curved, as they often contain tori.

\subparagraph{\bf 3.3.  Theorem.}
Let $\G=(G,\varphi)$ be a colored graph and let $n\geq 1$.  Then
$C_n(\G)$ satisfies the link condition.

\medskip 
We remark that H. Glover has also proved this theorem.

\begin{proof}
The fun in this proof lies in visualizing the link of
a cell; it is perhaps best illustrated by a picture.

Suppose the nine darkened cells of $\G$ in Figure \ref{npcthm.fig}
are given an ordering,
so that they represent a single (4-dimensional) cell $\r$ of $C_9(\G)$.
If the coloring allows us to include the edge
$e_x$ (shown dashed in the picture), we get a 5-cell 
$\s_x$ whose boundary contains $\r$; thus there is a vertex $x$
in the link of $\r$ representing the cell $\s_x$.  If the vertices $x$ and
$y$ are adjacent in the link of $\r$, then there is a 6-cell of $C_9(\G)$
whose boundary includes both $\s_x$ and $\s_y$; this happens exactly if
all the relevant cells are pairwise color-disjoint (and in particular
disjoint).

So, if $x$, $y$, and $z$ are vertices in the link of $\r$ which form
a triangle, then $e_x$, $e_y$, and $e_z$ (and the other darkened cells)
must be pairwise color-disjoint
in $\G$.  Therefore there is a 7-cell in $C_9(\G)$ formed (as a product)
from the darkened cells and $e_x$, $e_y$, and $e_z$.  This
7-cell is represented in the link of $\r$ by a 2-simplex spanned by $x$, $y$, 
and $z$, qed.

\bigskip
\begin{figure}[htb]\begin{center}
\psfrag{G}{\scalebox{2}{$\G$}}
\psfrag{a}{\scalebox{1.8}{$e_x$}}
\psfrag{b}{\scalebox{1.8}{$e_y$}}
\psfrag{c}{\scalebox{1.8}{$e_z$}}
\psfrag{d}{}
\psfrag{e}{}
\psfrag{f}{}
\scalebox{.5}{\includegraphics{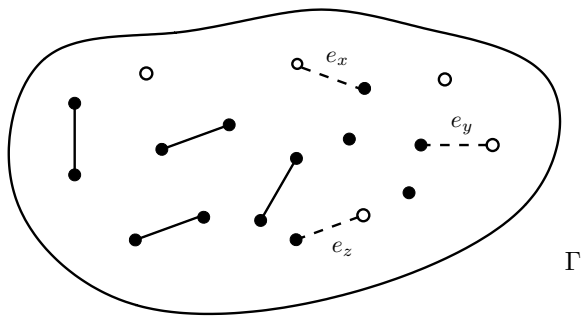}}
\caption{A colored graph $\G$ (the coloring is not shown).  The darkened cells
form a $4$-cell $\r$ in $C_9(\G)$.  If the coloring permits it, then adding 
a single dotted edge makes a $5$-cell, represented by a vertex in the link 
of $\r$.}
\label{npcthm.fig}
\end{center}
\end{figure}
\end{proof}

\subparagraph{\bf 3.4.  Corollary.}  The spaces $C_n(\G)$ are 
aspherical.\footnote{Strictly speaking, each connected
component of $C_n(\G)$
is aspherical.}  This is a property shared by all non-positively
curved spaces.

\section{Braid groups}

\subparagraph{\bf 4.1.}  The fundamental group of $C_n(\G)$, denoted
$P_n(\G)$, is called the \emph{pure ($n$-string) braid group} of $\G$.
For uncolored graphs $G$, it follows from {\bf 1.3} that if the 
edges of $G$ are subdivided enough, then $P_n(G)\iso \pi_1(C_n^{\rm top}(G))$
(see {\bf 1.2}).  For the purposes of this paper, we leave out the
choice of basepoint; in general, however, it matters.

In {\bf 1.6} and {\bf 1.8} we have $P_2(G)\iso\Z$.  In {\bf 2.3}
we have $P_2(\G)\iso \Z$ for any choice of basepoint.
In {\bf 1.9},
{\bf 1.10}, and {\bf 2.6}, $P_2$ is a surface group. 

\subparagraph{\bf 4.2.}
Fundamental groups of non-positively curved spaces have many 
famous properties.  For instance, they have solvable word problem
and conjugacy problem \cite{bridson}.  Since these spaces are aspherical
and
finite-dimensional, the braid groups of colored graphs
are also torsion-free.

\subparagraph{\bf 4.3.}
The symmetric group on $n$ letters acts freely on $C_n(\G)$ with
quotient denoted $UC_n(\G)$ (for \emph{unordered configurations}).
The fundamental group of $UC_n(\G)$ is the \emph{full $n$-string
braid group} $B_n(\G)$ of $\G$.  Paragraphs {\bf 2.5}, {\bf 3.3},
{\bf 3.4}, {\bf 4.2} hold with $C_n$ replaced everywhere by $UC_n$.

\subparagraph{\bf 4.4.}  As is well-known, the classical braid groups
are also fundamental groups of (aspherical) configuration spaces, this
time of points in the plane.
The braid groups are torsion-free and have solvable word and conjugacy
problems, as well as many other delightful properties.  Birman
\cite{birman} offers a good introduction.

\section{Graphs of groups}

\subparagraph{\bf 5.1.}  We now describe some additional structure
of the groups $P_n(\G)$.  Let $\G=(G,\varphi)$.
As in {\bf 1.9} let $p_1:C_n(\G)\to G$ be the projection 
onto the first coordinate.  For a cell $\s$ of $G$, define the colored
graph $\G_{\s}$ to be the colored subgraph of $\G$ spanned by those vertices
of $G$ which are colored differently from the vertices in $\bdry\s$
({\bf 1.4}).  Then for each vertex $v$ of $G$, $p_1^{-1}(v)$
is a copy of $C_{n-1}(\G_v)$, and
if $m$ is the midpoint (or any interior point) of an edge $e$, then 
$p_1^{-1}(m)$ is a copy of $C_{n-1}(\G_e)$.

Thus we can build $C_n(\G)$ inductively from $\{C_{n-1}(\G_v)
\,\big|\,v\ \mbox{a vertex of } G\}$ and $\{C_{n-1}(\G_e)\times I
\,\big|\,e\ \mbox{an edge of } G\}$ by using inclusion maps to glue 
each end of $C_{n-1}(\G_e)\times I$ to the appropriate $C_{n-1}(\G_v)$.

\subparagraph{\bf 5.2.  Example.}  $K_5$ has five vertices attached by
ten edges.  $C_2(K_5)$ is built from five copies
of $K_4$ attached with ten tubes $K_3\times I$.  

Cutting each of
the tubes gives a ``shirt decomposition'' of the surface $C_2(K_5)$.
Figure \ref{C2K5.fig}(b) shows one of the shirts.

\subparagraph{\bf 5.3.}  The key ingredient in establishing 
that {\bf 5.1} describes a splitting of $C_n(\G)$ as a 
\emph{graph of spaces} \cite{scottwall} over the graph $G$
is showing that each of the inclusion
maps $i:C_{n-1}(\G_e)\to C_{n-1}(\G_v)$ is 
$\pi_1$-injective.\footnote{Actually,
we mean that the inclusion is $\pi_1$-injective on each component
of $C_{n-1}(\G_e)$.}  This
follows from the geometric fact (which isn't hard to prove) that
$i(C_{n-1}(\G_e))$ is a \emph{locally convex} subset of $C_{n-1}(\G_v)$.

\subparagraph{\bf 5.4.}  We remark that this splitting is analogous
to the fiber bundle structure on the configuration spaces of points
in the plane \cite{birman}.  This bundle is the starting point 
for many arguments about braid groups.

\subparagraph{\bf 5.5.}  If all the spaces $C_{n-1}(\G_{\s})$ ({\bf 5.1})
are connected, then in fact $P_n(\G)$ splits as a graph of groups over $G$.
This is the case in {\bf 5.2}; the graph of groups diagram looks like
this ($F_3$ denotes a rank three free group):

\begin{figure}[h]\begin{center}
\psfrag{f}{\scalebox{1.6}{$F_3$}}
\psfrag{z}{\scalebox{1.6}{$\Z$}}
\scalebox{.6}{\includegraphics{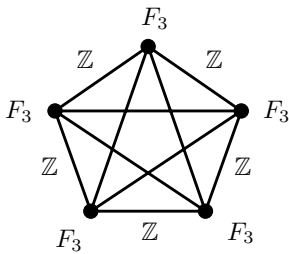}}
\caption{The surface group $P_2(K_5)$ splits as a graph of groups
as shown.  All edge groups, labelled or not, are $\Z$'s.}
\label{C2K5gog.fig}
\end{center}
\end{figure}

\begin{figure}
\begin{center}
\psfrag{A}{\scalebox{1.5}{$A$}}
\psfrag{B}{\scalebox{1.5}{$B$}}
\psfrag{C}{\scalebox{1.5}{$C$}}
\psfrag{O}{\scalebox{1.5}{$O$}}
\psfrag{C2}{\scalebox{1.5}{$C_2(Y)$}}
\psfrag{Y2}{\scalebox{1.5}{$Y^{(2)}$}}
\psfrag{Y}{\scalebox{1.5}{$Y$}}
\scalebox{.5}{\includegraphics{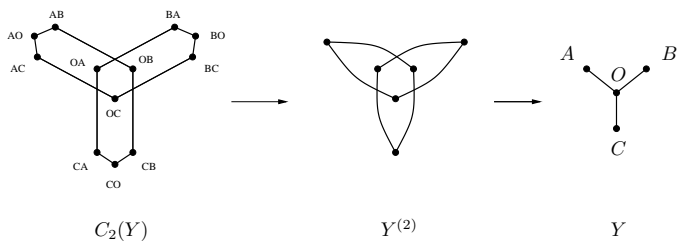}}
\caption{The space $C_2(Y)$ (left) splits as a graph
of \emph{connected} spaces over $Y^{(2)}$ (middle), which in turn maps
to $Y$ (right).  The composition of the two maps is $p_1$.
The group $P_2(Y)$ splits as a graph of groups
over $Y^{(2)}$, with all vertex and edge groups trivial.}
\label{C2K13gog.fig}
\end{center}
\end{figure}

\subparagraph{\bf 5.6.}  In general $P_n(\G)$ splits over a graph
$G^{(n)}$ comprised of a vertex for each component of each space
$C_{n-1}(\G_v)$ and an edge for each component of each space
$C_{n-1}(\G_e)$, attached in the natural way.  For example let $Y=K_{1,3}$.
The graph $Y^{(2)}$ is shown in Figure \ref{C2K13gog.fig} (middle),
and $P_2(Y)\iso\Z$ splits as a graph of groups over 
$Y^{(2)}$ with all vertex and edge groups trivial.  Compare with
{\bf 1.8}.

\bigskip

\section{More examples}\label{more}

\begin{figure}
\begin{center}
\scalebox{.7}{\includegraphics{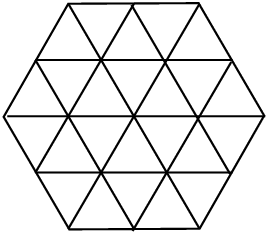}}
\caption{The link of a vertex of $C_3(K_7)$ is a torus.  Opposite edges of
the hexagon are identified.}
\label{toruslink}
\end{center}
\end{figure}

\subparagraph{\bf 6.1.  Example.}  Consider $C_3(K_7)$.  Much like 
$C_2(K_5)$
({\bf 1.9}), it is easy to check that each 2-cell lies in the 
boundary of exactly two 3-cells, and that the link of each 1-cell
is a (connected) 1-manifold.  Thus $C_3(K_7)$ is a 3-manifold away
from the vertices.  The link of each vertex is a torus, triangulated
as in Figure \ref{toruslink}.

If we remove the vertices of $C_3(K_7)$, we are left with an open
3-manifold $M$.  One can metrize each 3-cell as a regular ideal
cube in hyperbolic 3-space; this exhibits $M$ as a complete finite-volume
cusped hyperbolic 3-manifold.  $M$ contains many incompressible
surfaces, for instance the genus 6 surfaces $p_1^{-1}(x)$, where $x$ is
an interior point of an edge of $K_7$.

\subparagraph{\bf 6.2.  Example.}  The space $C_3(K_{4,4})$ is also
a 3-manifold except at some of its vertices.  This time some of the 
vertices are ``finite'' and have 2-sphere links, while others are
``ideal'' and have torus links as in Figure \ref{toruslink}.  Each 
3-cell of $C_3(K_{4,4})$ has two ideal and six finite vertices.  It
is an amusing exercise in 3-dimensional hyperbolic geometry to find a 
(``semi-ideal'') cube in hyperbolic space of the right shape to use
to metrize the cubes of $C_3(K_{4,4})$ so that $C_3(K_{4,4})$ minus
its ideal vertices is a complete finite-volume (cusped) hyperbolic
3-manifold.

\subparagraph{\bf 6.3.}  The next theorem generalizes Examples
{\bf 1.9}, {\bf 1.10}, {\bf 6.1}, and {\bf 6.2}.  It is proved
by basic combinatorial reasoning.  Recall ({\bf 2.6})
that the theorem doesn't hold for colored graphs. 

\medskip
{\bf Theorem.}   Let $G$ be a connected uncolored graph
and let $n\geq 2$.  The space $C_n(G)$ is homeomorphic to an
$n$-manifold away from the $(n-3)$-skeleton if and only if
$G=K_{2n+1}$ or $G=K_{n+1,n+1}$.

\subparagraph{\bf 6.4.  Corollary.}  Let $G$ be a connected uncolored graph.
Then $C_n(G)$ is homeomorphic to a closed $n$-manifold if and only if 
(i)  $n=1$ and $G$ is homeomorphic to a circle; or
(ii)  $n=2$ and $G=K_5$ or $G=K_{3,3}$.

\begin{proof}  By {\bf 6.3}, it suffices to check that for $n\geq 3$,
$C_n(K_{2n+1})$ and
$C_n(K_{n+1,n+1})$ contain $(n-3)$-cells $\s$ such that the link of $\s$
is a torus.  They do.
\end{proof}

\subparagraph{\bf 6.5.  Corollary.}  Let $G$ be an uncolored graph.
If $G$ contains $K_{2n+1}$ or $K_{n+1,n+1}$ as a subgraph, then
$C_n(G)$ has homological dimension $n$.

\begin{proof}  The $n$-cells of $C_n(K_{2n+1})$ or $C_n(K_{n+1,n+1})$
form a non-trivial $n$-cycle in the group of $n$-chains of
$C_n(G)$.
\end{proof}

\section{The outside world}

\subparagraph{\bf 7.1.}
Configuration spaces arise frequently within and without topology.
We have already mentioned the classical braid groups ({\bf 4.4}); much of our 
motivation comes from trying to compare and contrast
the configuration spaces (and braid groups) of the plane with the 
configuration spaces (and braid groups) of graphs.

Additionally, in 1933 van Kampen \cite{vankampen,shapiro}
generalized Kuratowski's criterion for
planarity of graphs by using the space $C_2(X)$ (where $X$ is an arbitrary
$n$-dimensional cell complex) to describe an obstruction to embedding
the $n$-complex $X$ into $\R^{2n}$.  In the literature $C_2(X)$ is
called the ``deleted product.''  We point out that if $X$ is not 
1-dimensional, then $C_2(X)$ may not be non-positively curved, even
if $X$ is.

\subparagraph{\bf 7.2.  Robotics.}  
It is probably not surprising that outside of topology there are many
situations in which an understanding of the possible motions of points
on a graph may be useful.  We briefly mention two examples.
 
In robotics, a common problem is this:  in a large factory, one has 
a collection of mobile robots which move along a network of tracks.
The problem is to design motion control algorithms which can move safely
from a given initial configuration of robots to a desired final
configuration.  Safely, of course, means that the robots shouldn't
collide.  Ghrist and Koditschek \cite{ghristkoditschek}
have had success solving this problem by using flows on configuration
spaces of graphs.  See also \cite{reconfig,monthly}.

\subparagraph{\bf 7.3.  Random walks.}
Second, the following problem of P. Winkler's \cite{tc98,winkler}
arose from considerations 
in theoretical computer science.  Consider two tokens taking discrete
walks on a fixed graph.  Suppose there is a demon whose goal
is to keep the tokens apart forever.  The demon decides, at each tick
of the clock, which token will move; moreover, \emph{the demon knows
exactly which (infinite) path each token will take!}  The problem
is to determine whether the ``clairvoyant demon'' has a positive
chance of succeeding if the walks are chosen randomly.

This problem can be formulated in terms of the 2-point configuration
space of the graph, and some results of \cite{winkler} and \cite{tc98}
can be duplicated in this way.  I hope
to use this approach to extend the known theorems.

\bibliography{../../../../Bibliography/Masterbib}
\bibliographystyle{amsplain}

\end{document}